\documentclass{amsart}
\usepackage{amscd,amssymb,amsxtra}

\numberwithin{equation}{section}

\theoremstyle{plain}
\newtheorem{Thm}{Theorem}
\newtheorem{Lem}[Thm]{Lemma}

\theoremstyle{definition}
\newtheorem{Def}[Thm]{Definition}

\newcommand{\NN}{\mathbb{N}}
\newcommand{\QQ}{\mathbb{Q}}

\newcommand{\calL}{\mathcal{L}}
\newcommand{\calG}{\mathcal{G}}
\newcommand{\calT}{\mathcal{T}}

\DeclareMathOperator*{\Wr}{Wr}

\begin{document}

\title{Hopfian Groups are Complete co-Analytic}
\author{Brian Pinsky}
\address
{Mathematics Department \\
Rutgers University \\
110 Frelinghuysen Road \\
Piscataway \\
New Jersey 08854-8019 \\
USA}
\email{bp466@math.rutgers.edu}

\begin{abstract}
Resolving a question posed in \cite{cpltSimon}, we show the set of Hopfian countable groups is a complete co-analytic subset of the standard borel space of countable groups.
\end{abstract}

\maketitle

\section{Introduction} \label{S:intro}
\begin{Def}
A group $G$ is \textbf{Hopfian} if every surjective group homomorphism from $G$ to $G$ is an isomorphism.
\end{Def}
Let $\calG$ be the standard Borel space of countably infinite groups, and let $\calG_{Hop}=\{G\in \calG | \text{$G$ is Hopfian}\}$. Then it is clear that $\calG_{Hop}$ is a co-analytic subset of $\calG$. 

\begin{Thm} \label{T:main}
$\calG_{Hop}$ is a $\Pi_1^1$ complete subset of $\calG$.
\end{Thm}

Let $\calT$ be the standard Borel space of trees on $\omega$.  To show $\calG_{Hop}$ is complete co-analytic, we will define a Borel map from $\calT$ to $\calG$ such that each tree in $\calT$ is well founded if and only if the corresponding group in $\calG$ is Hopfian.

As in \cite{cpltSimon}, we will start by mapping $\calT$ into $\calL$, the Borel space of linear orders of $\omega$, using the following result.

\begin{Thm}[Friedman-Stanly, \cite{fstrees}]
    There is a borel map $\varphi\colon \calT \to \calL$ such that, for each $T\in \calT$,
    \begin{itemize}
        \item $\varphi(T)$ is a pseudo-well ordering (that is, there is no infinite descending sequence hyperarithmetic in $\varphi(T)$)
        \item $\varphi(T)$ is a well ordering if and only if $T$ is well founded
    \end{itemize}
\end{Thm}

Next, we will define a Borel map from pseudo-well orders $W$ to groups $G_W$.  For this, we make use of a generalized wreath product construction from Silcock in \cite{silcock} to encode the order of $W$ into the lattice of subgroups of $G_W$.  

\section{Generalized wreath products} \label{S:wreath}

Let $\Lambda$ be a partial order, and let $\{H_\lambda | \lambda\in \Lambda\}$ be a $\Lambda$-indexed family of groups.  Let $S$ be the set $\left\{\left. x = (x_\lambda)\in \prod\limits_{\lambda \in \Lambda} H_\lambda \right| \text{$x_\lambda=1_{H_\lambda}$ for all but finitely many $\lambda$} \right\}$.  For each $\lambda\in \Lambda$ and $h\in H_\lambda$, we let $\xi_h$ be the permutation of $S$ defined by
$$
(\xi_h (x))_\mu =\begin{cases}
    x_\mu & \mu\neq \lambda\\
    x_\lambda & \mu = \lambda \text{ and } (\exists \eta>\lambda)\, x_\eta \neq 1_{H_\eta}\\
    h\cdot x_\lambda & \mu = \lambda \text{ and } (\forall \eta>\lambda)\, x_\eta =1_{H_\eta}
\end{cases}
$$
The wreath product $\Wr_{\lambda\in \Lambda}H_\lambda$ is the group of permutations of $S$ generated by all such $\xi_h$.  When $\Lambda=2$, $\Wr_{\lambda\in \Lambda}H_\lambda$ is the usual restricted wreath product $H_0\wr H_1$; and when no elements of $\Lambda$ are comparable, $\Wr_{\lambda\in \Lambda}H_\lambda$ is the direct sum of the $H_\lambda$.

If $\Gamma \subseteq \Lambda$, we will write $H_\Gamma$ for the subgroup of $\Wr_{\lambda\in\Lambda} H_\lambda$ generated by elements $\xi_h$ for $h\in H_\gamma$ for $\gamma\in \Gamma$.  Notice that the groups $H_\Gamma$ and $\Wr_{\gamma\in \Gamma} H_\gamma$ are isomorphic (where $\Gamma$ inherits its order from $\Lambda$). Note that $H_\Gamma$ is not a normal subgroup; for example, if $\Lambda=\{0,1\}$ and $\Gamma = \{0\}$, then $H_\Gamma$ will be one component of the base group in $H_0\wr H_1$.

In order to understand the normal subgroups of $H_\Lambda$, we will need to add the additional condition that $\Gamma$ is downwards closed in $\Lambda$.  To see why this condition needs to be imposed, let $\sim_\Gamma$ be the equivalence relation on $S$ defined by
$$
x\sim_\Gamma y \quad \Longleftrightarrow \quad (\forall \lambda\in \Lambda\setminus \Gamma)\, x_\lambda=y_\lambda.
$$ 
Since $\Gamma$ is downwards closed in $\Lambda$, the full group $H_\Lambda$ will 
permute the $\sim_\Gamma$ equivalence classes\footnote{
To see this, pick $x,y\in S$, some generator $h\in H_\Lambda$, and suppose $x\sim_\Gamma y$.  We want to show $hx\sim_\Gamma hy$.  Pick $\lambda\in \Lambda\setminus \Gamma$.  If $(hx)_\lambda\neq (x)_\lambda$, then $\forall \mu>\lambda. x_\mu=1$.  Since $\Lambda\setminus\Gamma$ is upwards closed and $x\sim_\Gamma y$, $\forall \mu>\lambda. y_\mu=1$, so $(hy)_\lambda = (hy)_\lambda$ since $(x)_\lambda=(y)_\lambda$.}.
Let $D_\Gamma$ be the kernel of the action of $H_\Lambda$ on $S/\sim_\Gamma$; i.e.\, the subgroup of $H_\Lambda$ containing all permutations of $S$ which fix each coordinate outside of $\Gamma$.  Then it is easily checked that $D_\Gamma$ is the normal closure of $H_\Gamma$ in $H_\Lambda$, and that the quotient $H_\Lambda/D_\Gamma\cong H_{\Lambda\setminus \Gamma}$.

In \cite{silcock}, Silcock shows that if every $H_\lambda$ is a non-abelian simple group, then every normal subgroup of $H_\Lambda$ has the form $D_\Gamma$ for some downwards closed $\Gamma\subseteq \Lambda$.  

\section{The proof of Theorem \ref{T:main}}

For each pseudo-well-order $W$, let $W^{op}$ be the opposite linear order; and for each $\lambda\in W^{op}$, let $H_\lambda$ be the alternating group $A_5$.  Let $G_W = H_{W^{op}} = \Wr_{\lambda\in W^{op}}H_\lambda$ be the wreath product, defined in Section \ref{S:wreath}.  Then $G_W$ is clearly a countable group, and the mapping $W\mapsto G_W$ is Borel.

\begin{Lem}
    If $W$ is a well-order, then $G_W$ is Hopfian
\end{Lem}
\begin{proof}
    Since $W$ is a well-order, each non-empty downwards closed subset of $W^{op}$ has a unique maximum element; so the set of non-trivial normal subgroups of $G_W$ has order type $W^{op}$.
    
    Suppose that $f\colon G_W\to G_W$ be a surjective group homomorphism. Then, by \cite{silcock}, $\ker(f)$ must be  $D_{X'}$ for some downwards closed $X'\subseteq W^{op}$.  Let $X=W^{op}\setminus X'$.  Then $X$ is an upwards closed subset of $W^{op}$; i.e.\ an initial segment of the well-order $W$ equipped with the opposite order.  Note that $G_W/\ker(f)\cong H_{X^{op}}$; and hence
    $$
    G_W \cong G_W/\ker(f)\cong H_{X^{op}} \cong G_X.
    $$
    
    It follows that the sets of normal subgroups of $G_W$ and $G_X$ have the same order type; so the well-orders $W$ and $X$ are isomorphic.  Since $X$ is an initial segment of $W$, it follows that $W=X$, and $\ker(f)$ is trivial. Thus $f$ is injective.
\end{proof}

For the converse, we will need a theorem of Harrison.

\begin{Thm}[Harrison, \cite{harrison}]
    If $R$ is a pseudo-well-order of $\NN$, and $R$ is not a well-order, then $R$ has order type $(\omega_1^R\times (1+\QQ))+\alpha$ for some $\alpha<\omega_1^R$
\end{Thm}
In particular, it follows that
every ill founded pseudo-well-order is isomorphic to a proper initial segment, since
$$
(\omega_1^R\times (1+\QQ))+\alpha \cong 
(\omega_1^R \times (1+(-\infty,0)\cap \QQ)) \sqcup (\alpha \times \{\, 0\,\}).
$$

\begin{Lem}
If the countable pseudo-well-order $W$ is not a well order, than $G_W$ is not Hopfian
\end{Lem}
\begin{proof}
    Let $L$ be a proper initial segment of $W$ such that $W\cong L$ and let $N=D_{W\setminus L}$ be the normal subgroup of $G_W$ corresponding to initial segment $W\setminus L$ of $W^{op}$.  Then $G_W/N$ is isomorphic to $G_L$; and since $L$ and $W$ are isomorphic, it follows that $G_L \cong G_W$.  Thus $G_W$ is not Hopfian.
\end{proof}

This completes the proof of Theorem \ref{T:main}.

\section{Further Questions}
\begin{Def}
    A group $G$ is co-Hopfian if every injective group homomorphism from $G$ to $G$ is an isomorphism
\end{Def}

Once again, it is clear that the class of co-Hopfian groups is also a co-analytic subset of $\calG$, and it was conjectured in \cite{cpltSimon} to be complete co-analytic.  Since every countable linear order is not co-Hopfian (see \cite{noCoHopOrders}), proving this will require more than simply encoding linear orders in groups.  Currently, we do not know any family of structures for which the co-Hopfian ones are $\Pi^1_1$ complete. However, we do know that if such a family exists, then the co-Hopfian groups should also be $\Pi_1^1$ complete, using an encoding similar to that in \cite{cpltSimon}.


\begin{thebibliography}{99}

\bibitem{noCoHopOrders}
Ben Dushnik and E. W. Miller
{\em Concerning similarity transformations of linearly ordered sets\/}, 
Journal of Amer. Math. Soc. {\bf 46} (1940), 322--326

\bibitem{fstrees}
H. Friedman and L. Stanley, 
{\em A Borel reducibility theory for classes of countable structures\/},
Journal of Symbolic Logic {\bf 54} (1989), 894--914.

\bibitem{harrison}
J. Harrison, 
{\em Recursive pseudo-well-orderings\/}, 
Trans. Amer. Math. Soc. {\bf 131} (1968), 526--543. 

\bibitem{silcock}
Howard L. Silcock,
{\em Generalized wreath products and the lattice of normal subgroups of a group\/}, 
Algebra Universalis {\bf 7} (1977), 361--372

\bibitem{cpltSimon}
Simon Thomas,
{\em Complete groups are complete co-analytic\/},
Archive for Mathematical Logic {\bf 57} (2018), 601--606

\end{thebibliography}
\end{document}